\documentclass[1pt,notitlepage,twoside,a4paper]{amsart}

\usepackage{amsmath,amssymb,enumerate}

\usepackage{epsfig,fancyhdr,color}

\usepackage{amssymb}
\usepackage{amsmath,amsthm}    
\usepackage{latexsym}
\usepackage{amscd} 
\usepackage{psfrag}
\usepackage{graphicx} 
\usepackage[latin1]{inputenc}  
\usepackage[all]{xy} 
\usepackage[mathcal]{eucal}

\theoremstyle{definition}

\theoremstyle{remark}

\def\interieur#1{\mathord{\mathop{\kern 0pt #1}\limits^\circ}}

\definecolor{NoteColor}{rgb}{1,0,0}

\title[Mathematics, physics and philosophy
in Riemann's work]{Mathematics, physics and philosophy
in Riemann's work and beyond}

\author{Athanase Papadopoulos}
\address{A. Papadopoulos, Institut de Recherche Math{\'e}matique Avanc\'ee,
Universit{\'e} de Strasbourg and CNRS,
7 rue Ren\'e Descartes,
 67084 Strasbourg Cedex, France and Brown University, Mathematics Department, 
 151 Thayer Street
Providence, RI 02912, USA.}

 \date{\today}

\begin{document}

\begin{abstract}
This is the introduction I wrote for the multi-authored book \emph{From Riemann to differential geometry and relativity}, L. Ji, A. Papadopoulos and S. Yamada, (ed.), Berlin, Springer, 2017. 
The book consists of twenty chapters, written by various authors. This introduction, besides giving the information on the content of the book, is a quick review of the topics on which Riemann worked and of the impact of this work on mathematics (topology, complex geometry, algebraic geometry, integration, trigonometric series, Riemannian geometry, 	etc.), philosophy and physics. 
 
\bigskip


\end{abstract}
  \maketitle

              The present book is an addition to the living literature on Riemann. It contains a series of introductory essays in which the authors comment on some of Riemann's writings with the goal of making them more accessible, followed by surveys of some recent research topics rooted in Riemann's work or strongly motivated by his ideas. The overall goal is to give a comprehensive overview of Riemann's work, the origin of his ideas and their impact on mathematics, philosophy and physics. The various authors -- each one with his own style -- get into a great variety of subjects including Riemann surfaces, elliptic and Abelian integrals, the hypergeometric series, differential geometry, topology, integration theory, the zeta function, minimal surfaces, uniformization, trigonometric series, electromagnetism, heat propagation, Riemannian Brownian motion, and several other topics to which Riemann made essential contributions or that were greatly influenced by his work. 
              One difference between this book and the existing books on Riemann is that it contains a significant part devoted to Riemann's impact on philosophy  (there are three chapters on this subject out of a total of twenty chapters in the book) , while another consequential part (again, three chapters) is concerned with the impact of Riemann's ideas on the theory of relativity. Let us add that even though part of  the book deals with subjects that are treated in other books on Riemann, it is always useful to have, in the mathematical literature, surveys on the same subject written by different persons, each survey reflecting its author's interests and his ideas on what is important and what is only secondary material (although in Riemann's case secondary material is very rare).
              
              Riemann's influential habilitation lecture,\index{Riemann! habilitation lecture}\index{Habilitation lecture!Riemann} \emph{\"Uber die Hypothesen, welche der Geometrie zu Grunde liegen} (On the hypotheses that lie at the bases of geometry) is at the center of the discussion in several chapters of the book. The repercussion of this lecture in mathematics, physics and philosophy is immense. Occurrences of a single text that had such a profound influence on these three branches of human knowledge are very rare in history. Other examples may be found in the writings of very few thinkers: Aristotle, Newton, Leibniz,\index{Leibniz, Gottfried Wilhelm (1646--1716)}  Descartes and Poincaré are some of them, and it is difficult to find more names.

    The present introduction has several purposes. The first one is to provide the reader with a short summary of the topics that he will find in this book. Reading this summary will  give him an idea on the great variety of themes on which Riemann worked and on the impact of his ideas. Another purpose, on which we stress in the last part of the introduction, is to transmit a few thoughts, beyond Riemann's work, on the intricate relation between mathematics, physics and philosophy.    
  
  \medskip
     
 The volume is organized into a preamble, four main parts, and a concluding chapter.  Each of the four parts contains a series of essays, arranged in chapters. 
  
  \medskip
 
 Chapter 1, which constitutes the preamble, written by the author of the present introduction, is  an overview of the prehistory of some of the main mathematical fields on which Riemann worked. In other words, the chapter concerns the predecessors of Riemann, more precisely, the mathematicians who started the fields in which he worked, and those who exerted a major influence on him. It appears clearly from this overview that for most of the questions which Riemann addressed, Euler stands in the background, as a dominating figure. This concerns the theory of functions (in particular algebraic functions and functions of a complex variable), of elliptic integrals, of Abelian integrals, of the hypergeometric series, of the zeta function and of Riemann's ideas on space, as well as his work on topology,  differential geometry, trigonometric series, and integration, and his use of the techniques of the calculus of variations. Even though Euler was not the initiator of all these fields, he is, for most of them, the one who brought them to matureness. This applies in particular to the theories of algebraic and complex functions, to that of elliptic integrals and to the zeta function. Despite the fact that the history of topology can be traced back to the Greeks, and then to Leibniz and Descartes, Euler was the first to solve topological problems with the conviction that these problems are proper to this field, and that the classical method of analysis and algebra are insufficient for their solution.  Another major figure to whom Riemann has an enormous debt for what concerns his intellectual and mathematical development, is his mentor Gauss, who worked on all the topics that Riemann tackled. We also know, from Riemann's writings and his correspondence, that he was a dedicated reader of Euler and Gauss's works. What we said about mathematics in Riemann's writings also holds for physics and philosophy, that is, it is possible to trace back several important ideas of Riemann in these fields to Euler and to Gauss.
 This chapter is also in some sense an essay on historical progression in mathematics and it is an occasion of revisiting the texts written by several pre-eminent mathematicians of which we are the heirs, and on whose shoulders we stand.

 Part I, composed of Chapters 2 to 8, is an exploration of Riemann's works and their impact on mathematics and physics. Some of these chapters have a historical character, and others contain detailed reviews of some of Riemann's published works. Some relations with the works of Riemann's contemporaries are also highlighted.

Part II,  containing Chapters 9 to 11, is more directed towards the philosophical aspects of Riemann's work. It focuses  in particular on his ideas on space, making relations with conceptions of other thinkers on the same subject, and exploring the impact of these ideas on modern philosophy. The notion of \emph{Mannigfaltigkeit}\index{Mannigfaltigkeit} (usually translated as multiplicity,\index{multiplicity} or manifoldness),\index{manifoldness} which existed in the philosophical language, which Riemann introduced in mathematics, and which is an ancestor of the modern notion of manifold, is thoroughly discussed.

Part III,  consisting of Chapters 12 to 16, is a series of five surveys of modern mathematical research topics that are based on ideas originating in Riemann's work. These  topics belong to the fields of conformal geometry, algebraic geometry, the foundations of mathematics, integration, and probability theory.

Part IV, consisting of Chapters 17 to 19, is a collection of surveys on the theory of relativity and more especially on questions in relativity that are directly inspired or that rely on Riemann's work.

 Chapter 20, the concluding chapter in this volume, is written by Lizhen Ji and is meant to give a quick overview of the life and works of Riemann. It contains in particular a brief summary of each of Riemann's published articles, together with a list of notions that Riemann introduced and that are named after him.  
\medskip

We now present in more detail the content of Parts I to IV, chapter by chapter.

\part*{Part I}

 Chapter 2 is written by Jeremy Gray and it starts by a review of some important aspects of Riemann's habilitation lecture,\index{Riemann! habilitation lecture}\index{Habilitation lecture!Riemann} \emph{On the hypotheses that lie at the bases of geometry}, from the points of view of mathematics, physics and philosophy, highlighting the consequences of Riemann's conception of space. It is a matter of fact that the philosophical notion of space has also implications on Riemann's mathematical works. According to Gray, Riemann's conception of space is related to the question of whether objects of geometry are described by coordinates or not. Following the line of work started by Gauss on the differential geometry of surfaces, Riemann formulated in a novel way the question of ``determination of position" in a manifoldness,\index{manifoldness} and that of the difference between intrinsic and extrinsic properties. 
 
 Gray then turns to the more general question of the interaction between the mathematical ideas expressed in Riemann's lecture and physics and philosophy. Riemann claims in his memoir that he was influenced by Herbart,\index{Herbart, Johann Friedrich (1776--1841)} without being explicit on that. Based on passages from Riemann's notes on philosophy, Gray presents some ideas of Herbart's\index{Herbart, Johann Friedrich (1776--1841)} on space, time, and motion, and he discusses  the way they were received by Riemann. For instance, Herbart\index{Herbart, Johann Friedrich (1776--1841)} addressed the question of whether or not our knowledge of space, time, and motion is generated by our experience. Gray explains where Riemann agrees or disagrees with Herbart's ideas. Like Newton before him, Riemann disputed the idea of an action at a distance. He imagined, like Euler did before him, that space is filled with a substance -- ether -- whose properties are responsible for the transmission of the forces of nature. This philosophico-physical idea was the driving force that led Riemann to the discovery of what became known as Riemannian geometry. From the physical viewpoint, this geometry is seen as the study of spaces with infinitesimal physical forces that are responsible for curvature. 
 
 Gray then considers another subject, namely, Riemann's theory of electrodynamics,\index{electrodynamics} formulated in his article  \emph{Ein Beitrag zur
Elektrodynamik} (A contribution to electrodynamics). The paper was presented to the Royal Academy of Sciences at G\"ottingen in 1858, but subsequently withdrawn, and it was published posthumously in 1867. This article by Riemann is also analyzed in Chapter 3 of the present volume. It is motivated by a question concerning the velocity of electrical interaction. Riemann argues that this velocity, ``within the limits of errors of observation, is equal to that of light." 

After electrodynamics, Gray comments on Riemann's theory of heat diffusion,\index{heat diffusion} expressed in another essay known as the \emph{Commentatio}, whose aim is to find conditions on the distribution of heat in an infinite, homogeneous, solid body under which a system of curves remains isothermal for an indefinite period of time. Riemann formulated this problem in terms of a positive definite quadratic form with constant coefficients at each point on the solid body, governing the heat flow.  He then interprets the physical problem mathematically, as a problem concerning the reduction of a quadratic form.  Thus, Riemann again places a physical problem at the heart of Riemannian geometry.  In the same chapter, Gray addresses the question of the influence of the \emph{Commenatio} on the works of later mathematicians.

Chapter 3, written by Hubert Goenner, deals with Riemann's work on electromagnetism. It is built around Riemann's paper  \emph{Ein Beitrag zur
Elektrodynamik}\index{electrodynamics} (A contribution to electrodynamics) which is also considered in Chapter 2.
Goenner analyzes this paper in some detail. He reminds the reader that the idea that electrical interaction is not instantaneous was voiced by Gauss already in 1845. Based on it, Riemann deduced an explanation of the electrodynamic actions of galvanic currents. Goenner highlights several important points in Riemann's paper, explaining how Riemann's theory anticipates that of Maxwell. He also mentions connections with later discoveries of Riemann that led him to change some important points in his theory. Riemann later on addressed these questions in more detail, in a course entitled ``The mathematical theory of gravitation, electricity and magnetism," that
 he gave in the summer 1861. 
 
 Goenner's commentary is a useful reading guide to Riemann's paper. It presents Riemann's work in a large perspective comprising the works of Gauss, Weber and others. Incidentally, Goenner provides an answer to the question of why Riemann withdrew his paper, namely, a wrong factor that Riemann included in a function under an integral sign. According to Goenner, Riemann realized, after submitting this paper, that a different factor should be there and this led him to withdraw the paper. There are other conjectural reasons for Riemann's withdrawal of the paper,  for instance, the fact that Riemann realized that he used a trivially forbidden interchange of integration -- this explanation is the one given by Gray in Chapter 2 of the present volume.

 Chapter 4, by Christian Houzel, concerns Riemann's solution of Jacobi's problem of inversion of Abelian integrals. These are integrals of the form $\int_{z_{0}}^z R(w,z)dz$ where $R(w,z)$ is a rational function of the two variables $w$ and $z$ that are related by an algebraic equation $f(w,z)=0$. The Jacobi inversion problem, which was formulated by Jacobi in 1832, generalizes the inversion problem for elliptic integrals to which Riemann also contributed in an essential way. The inversion of elliptic integrals leads to the so-called doubly periodic functions, that is, holomorphic functions defined on the torus. The inversion of Abelian integrals leads to what became known later on as automorphic functions, on more general surfaces. Elliptic integrals are in some sense generalizations of inverse trigonometric functions ($\int_0^x\frac{dt}{\sqrt{1-t^2}}$ represents the arcsine function, a special case of the class of elliptic integrals of the form $\int \frac{dx}{\sqrt{1-x^n}}$) and a major idea behind this study is that inverses of elliptic integrals may behave in some sense like trigonometric functions, having periods, addition formulae, etc. Weierstrass also worked on the Jacobi inversion problem. Riemann sketched a solution of this problem in his famous memoir  {\it Theorie der Abel'schen Functionen} (Theory of Abelian functions), written in 1857, without giving a complete proof. He completed the proof in his memoir {\it \"Uber das Verschwinden der $\vartheta$-Functionen} (On the vanishing of $vartheta$ functions), published 1866. Houzel makes a historical survey of this inversion problem and gives an outline of Riemann's proof. This proof uses all the concepts that Riemann introduced, including the representation of algebraic functions by Riemann surfaces that are coverings of the Riemann sphere, his formulation of the problem in terms of periods of differential forms of the first kind on the associated Riemann surface, and his use of what became known later on as Riemann's theta functions.  His 1857 memoir concludes with a proof of the fact that integrals of the algebraic differential forms on a Riemann surface may be expressed as quotients of products of translated theta functions. Riemann also contributed to the classification and the study of moduli of Abelian integrals.
 In the last section, Houzel indicates some later developments of Riemann's results by A. Weil (1948) and G. Kempe (1971--73).

 Chapter 5 by Sumio Yamada concerns Riemann's work on minimal surfaces. It consists of an overview, from a modern viewpoint, of Riemann's two papers  on the subject, \emph{\"Uber die Fl\"ache vom kleinsten Inhalt bei gegebener Begrenzung} (On the surface of least area with a given boundary) and \emph{Beispiele von Fl\"achen kleinsten Inhalts bei gegebener Begrenzung} (Examples of surfaces of least area with a given boundary). Both papers were finalized after Riemann's death by K. Hattendorff to whom Riemann had left a set of notes on the subject. At the same time, Yamada makes a comparison between Riemann's work and that of Weierstrass on the same subject. He shows that Riemann's notes contain several  
results on minimal surfaces which are now classical, including the Weierstrass-Enneper 
representation,  Schwarz's explicit construction of minimal surfaces, as well as the Schwarz--Christoffel transformation. He also mentions relations with the works of Euler and Lagrange and with Riemann's own work on the Riemann mapping theorem.\index{Riemann mapping theorem}\index{theorem!Riemann mapping}

 Chapter 6, by the author of this introduction, is a survey of the ideas from physics that are contained in Riemann's mathematical papers, and on the  mathematical problems that he tackled that were motivated by physics. In fact, it is sometimes not easy to separate Riemann's mathematical ideas from physics, and it is clear that for certain topics, Riemann did not make any difference between mathematics and physics. Furthermore, Riemann's philosophical ideas are often in the background of his work on mathematics and physics. The main goal of Chapter 6 is to try to convey this general theoe, by analyzing several writings of Riemann. These include his habilitation lecture, his habilitation text on trigonometric series,\index{Riemann! habilitation text}\index{Habilitation text!Riemann} the \emph{Commentatio}, a paper on differential geometry motivated by the problem of expressing the temperature at a point  of a homogeneous solid body in terms of time and a system of coordinate on the body,  his paper on the equilibrium of electricity, his paper on the propagation of planar air waves, his paper on the functions representable by Gauss's hypergeometric  series $F(\alpha,\beta,\gamma,x)$, and a few others. The chapter also contains a discussion of some of Riemann's philosophical ideas, mentioning several of Riemann's predecessors in this domain, in particular the Greek philosophers. The intricate relation between physics and mathematics in Riemann's work that is surveyed in this chapter  is a vast field.

 Chapter 7, also written by the author of the present introduction, is an essay on the works of Cauchy and Puiseux, the two French predecessors of Riemann on the theory of functions of a complex variable. 
 
  Cauchy started working in this field in 1814, that is, twelve years before Riemann was born. He introduced several concepts which were useful to Riemann, including line integrals, the dependence of such an integral on the 
homotopy class of the path of integration, and the calculus of residues. Cauchy wrote a large number of papers on this subject. In 1851, he discovered, for a complex function of a complex variable, the notion of derivative independent of direction, and he showed that the real and imaginary parts of such a function must satisfy the two partial differential equations that became known as the Cauchy--Riemann equations. At the end of the same year, Riemann submitted his doctoral dissertation, which contains the same concept of derivative, with the same characterization.

Puiseux was much younger than Cauchy and he followed his lectures on the theory of functions of a complex variable. He wrote two remarkable papers on this subject, and in particular on the question of uniformizing (that is, making single-valued) a multi-valued function defined implicitly by an algebraic equation. Puiseux' first paper was published in 1850, that is, one year before Riemann defended his doctoral dissertation in which he introduced the concept of Riemann surface. Interpreted in the right setting, this work of Puiseux inaugurates a group-theoretic point of view on the theory of Riemann surfaces. The second paper by Puiseux was published in 1851. The notions that Puiseux discovered constitute a combinatorial version of Riemann surfaces. Hermite made the relation between uniformization and Galois theory, based on the work of Puiseux. This is also reported on in Chapter 7.

Chapter 7 and the next one are also an occasion for the reader to learn about the lives of several pre-eminent mathematicians who florished at the epoch of Riemann and who had ideas close to his.
 Gaining insight into the life of a great mathematician is interesting even if this life has nothing exceptional. It often makes us understand his motivations and makes his work more familiar.

    Chapter 8, again by the author of the present introduction, is  concerned with the reception of the concept of Riemann surface by the French school and how this concept is presented in the French  treatises on analysis published in the few decades that followed Riemann's work on this subject.

It took several years to the mathematical community to understand the concept of Riemann surface that was conceived by Riemann as the base ground for general meromorphic functions and on which a multi-valued function becomes uniform, and to accept the validity of some major results that Riemann proved regarding these surfaces -- his theorem saying that a meromorphic function is determined by its singularities, and other results in the same vein. 

 At the same time, Chapter 8 is a survey of the remarkable French school of analysis that started with Lagrange, then Cauchy, and attained a high degree of maturity in the second half of the nineteenth century. We also comment on the relations between this school and the German one.
 The reader will find in this chapter a survey of works related to the concept of Riemann surface and related matters (elliptic and Abelian integrals, the topology of surfaces, the uniformization of multi-valued functions, etc.)  by various French authors including Briot, Bouquet, Appell, Goursat, Picard, Simart, Fatou, Jordan, Halphen, Tannery, Molk, Lacour and Hermite.
  
  \part*{Part II}

 The question of space,\index{space} which was already addressed several times in  Part I, is thoroughly studied in the three chapters that constitute Part II of the present volume.
  
 Elaborate psychological, physical, philosophical and mathematical theories of space\index{space} were developed by various thinkers since the birth of philosophy. We live in a space, even in the few months before we are officially born.\index{space} Everyone has a feeling of space, and we are supposed to have the impression that this space is Euclidean. Making a philosophy of space implies going a step further than these primary feelings. The Greeks, since the early Pythagoreans, wondered about the properties of space that go beyond our immediate senses, addressing the questions of whether space is infinite, whether it is full of matter, whether void exists, etc.  Riemann had his own ideas on space, and these are contained in his habilitation lecture,\index{Riemann! habilitation lecture}\index{Habilitation lecture!Riemann} and also in unfinished notes that were published pusthumously. The mathematical notion of manifold was born from Riemann's reflections on space.\index{space}  This is the major theme addressed in Chapters 9 to 11.
    
   Chapter 9, by Ken'ichi Ohshika, is essentially concerned with the notion of manifold, starting from the first introduction by Riemann, in his habilitation lecture,\index{Riemann! habilitation lecture}\index{Habilitation lecture!Riemann} in a mathematical context (but still with a high philosophical flavor), of the word \emph{Mannigfaltigkeit},\index{Mannigfaltigkeit} usually translated by \emph{multiplicity}\index{multiplicity}. The survey takes us until the modern notion of manifold, developed in the twentieth century, including the introduction of the specialized notions of Hausdorffness and differentiability. The contributions of Hilbert, Weyl, Kneser, Veblen--Whitehead and Whitney are surveyed. Poincar\'e's two definitions of a manifold, formulated at the turn of the nineteenth century, are also presented. The first definition is close to that of a submanifold, and the other one, using the notion of analytic continuation, is closer to the modern definition of a manifold using charts. The philosophical background of Riemann is also discussed, including the influences of Kant and Herbart\index{Herbart, Johann Friedrich (1776--1841)} on his ideas. Ohshika explains how Riemann's point of view differs from that of Kant, who regarded Euclidean space and its geometry as given a priori, thus excluding in principle the concepts of non-Euclidean geometries, and who apparently never thought of a possibility of alternative views on space and time.

  Chapter 10 by Franck Jedrzejewski is mainly devoted to the influence of Riemann on two pre-eminent twentieth-century French thinkers, Gilles Deleuze\index{Deleuze, Gilles (1925--1995)} and Félix Guattari.\index{Guattari, F\'elix (1930--1992)} 
 
 Deleuze\index{Deleuze, Gilles (1925--1995)} was a philosopher with a large spectrum of themes, including literature, politics, psychoanalysis\index{psychoanalysis} and art. Guattari\index{Guattari, F\'elix (1930--1992)} was a philosopher and a psychoanalyst\index{psychoanalysis} who followed during several years the famous seminar led by Lacan, who at the same time was Guattari's psychoanalyst.\index{psychoanalysis}  Deleuze and Guattari\index{Guattari, F\'elix (1930--1992)}  had a long and fruitful collaboration which culminated in their joint work \emph{Capitalism and Schizophrenia}, a complex philosophical essay in two parts entitled \emph{Anti-Oedipus}\footnote{G. Deleuze\index{Deleuze, Gilles (1925--1995)} and F. Guattari,\index{Guattari, F\'elix (1930--1992)}   L'anti-\OE dipe, Paris, \'Ed. de Minuit, 1972. English translation by  R. Hurley, M. Seem and H. R. Lane: Anti-\OE dipus, London and New York: Continuum, 2004.} (1972) and \emph{A Thousand Plateaus}\footnote{G. Deleuze\index{Deleuze, Gilles (1925--1995)} and F. Guattari,\index{Guattari, F\'elix (1930--1992)}   Mille Plateaux, Paris, \'Ed. de Minuit, 1980.  English translation by B. Massumi: A Thousand Plateaus, London and New York: Continuum, 2004.} (1980).   In this work,  the authors address various questions concerning political action, desire, psychology, economics, society, history and culture. As the name of the first volume suggests, the work is critical of psychoanalysis\index{psychoanalysis} as it was conceived by Freud. In fact, throughout his relation with Deleuze,\index{Deleuze, Gilles (1925--1995)} Guattari\index{Guattari, F\'elix (1930--1992)} distanced himself from Lacan,\index{Lacan, Jacques (1901--1981)} and Deleuze\index{Deleuze, Gilles (1925--1995)} and Guattari\index{Guattari, F\'elix (1930--1992)} expressed their disagreement with the fact of reducing the unconscious mind to the family circle of the individual (his relation with his parents).  The second volume of  \emph{Capitalism and Schizophrenia} contains a discussion, evaluation and critique of works of Freud, Jung,  Reich and Francis Scott Fitzgerald.  The publication of the two volumes generated a large debate in the intellectual milieu in France and sometimes beyond, and the ideas formulated by  Deleuze\index{Deleuze, Gilles (1925--1995)} and Guattari\index{Guattari, F\'elix (1930--1992)} had a non-negligible political influence in the last quarter of the twentieth century. Their work belongs to the so-called post-structuralist and  transcendental empiricism postmodernist currents.

At this point, the reader may rightly ask: \emph{What does all this have to do with Riemann and with mathematics?}
Another question will also soon be addressed:  \emph{Why were some twentieth-century French philosophers interested in Riemann and how were they influenced by him?}
 In Chapter 10  Jedrzejewski brings some answers to these questions. As a preliminary attitude, the reader has to realize that in the same way as there are mathematicians interested in philosophy, there are philosophers interested in mathematics, and this has been so since antiquity. Not only these philosophers were interested in mathematics, but they brought mathematical notions and ideas into the realm of philosophy, and they used them in their works, sometimes as essential elements in formulating systems of thought which they wanted to be coherent and built on a logical basis. We can quote here Jules Villemin (1920--2001), another pre-eminent French philosopher, from his major work \emph{La Philosophie de l'algèbre}:\footnote{Presses universitaires de France, 1962.}
 \begin{quote}\small
 There exists an intimate -- although less apparent and more uncertain -- relationship between pure mathematics and theoretical philosophy. History of mathematics and of philosophy shows that a renewal of the methods of one of them, each time had an impact on the other one.\footnote{Il existe un rapport intime quoique moins apparent et plus incertain entre les Mathématiques pures et la Philosophie théorique. L'histoire des mathématiques et de la philosophie montre qu'un renouvellement des méthodes de celles-là a, chaque fois, des répercussions sur celles-ci.}
 \end{quote}
 
 There are many instances in the history of ideas of works on philosophy having an impact on the development of mathematics. We recall for example Plato's influence on the development of geometry and Aristotle's influence on axiomatics and the foundations of mathematics.  The various views from which philosophers considered the notion of space had also a certain impact on mathematics, and this theme is considered in the various chapters that constitute the  second part of the present volume. One may also mention the enormous influence of such views on the research conducted during several centuries on Euclid's parallel axiom that led eventually to the discovery of non-Euclidean geometry. 
 
 Deleuze's\index{Deleuze, Gilles (1925--1995)} philosophical theories are rooted in the works of mathematicians like Riemann, Leibniz, Whitehead, Albert Lautman and Gilles Ch\^atelet.  Already in 1968, in an essay entitled \emph{Diff\'{e}rence et r\'{e}p\'{e}tition}, he expressed the fact that an idea, from the point of view of its organization, is the philosophical analogue of a continuous multiplicity\index{multiplicity} in the sense of Riemann. Guattari developed a philosophical concept of multiplicity,\index{multiplicity} based on Riemann's \emph{Mannigfaltigkeit},\index{Mannigfaltigkeit} as an alternative to the notion of \emph{substance},\index{substance} which is one of the key concepts in metaphysics. Many other mathematical terms, like dimension, continuity, variability, order, and metric, acquired a philosophical significance in Deleuze's\index{Deleuze, Gilles (1925--1995)} work.

In Chapter 10, Jedrzejewski makes a detailed comparison between some texts of Riemann and those of Deleuze\index{Deleuze, Gilles (1925--1995)} and Guattari.\index{Guattari, F\'elix (1930--1992)} In fact, many twentieth-century philosophers addressed questions related to or arising from mathematics, its logic and its language. Deleuze\index{Deleuze, Gilles (1925--1995)} was particularly fascinated by topology.  He was influenced of Leibniz, relying at the same time on his metaphysics, his differential calculus and his ideas on topology. Jedrzejewski also mentions the work of another French philosopher, Henri Bergson.\footnote{The French philosopher and teacher Henri Bergson (1859--1941) was awarded in 1927 the Nobel prize for literature. He had a mathematical background and there is a famous controversy between him and Einstein\index{Einstein, Albert (1879--1955)} concerning the philosophical notion of time, which might be interesting for the reader of this book.}

Beyond its relation with Riemann, this chapter by Jedrzejewski is an interesting example of how mathematics meets philosophy. The chapter is written in French, the reason being that Jedrzejewski wanted to use the original Deleuzian terms. An extended English summary of the chapter is provided by the author.

Chapter 11,  by Arkady Plotnistky, has also a philosophical character. It concerns the ``conceptual" nature of Riemann's thinking and its implications in mathematics, physics and philosophy. The word ``concept" is used here in a technical sense explained by Plotnitsky, who relies on another philosophical essay by Deleuze\index{Deleuze, Gilles (1925--1995)} and Guattari,\index{Guattari, F\'elix (1930--1992)} \emph{What is philosophy?} (1994), in which these authors view thought (``la pensée"), with its creative nature,  as a confrontation between the brain and chaos.   
Plotnitsky's discourse is at the level of ``concepts of concept,"  promoted by Deleuze\index{Deleuze, Gilles (1925--1995)} and Guattari\index{Guattari, F\'elix (1930--1992)} in the realm of philosophical thinking, transferred (by Plotnitsky) to the physical and mathematical worlds as well, despite the fact that these authors claim that their concept of concept pertains uniquely to philosophical thinking. The discussion around this concept of concept\index{concept} and the confrontation between the ideas of Riemann, Hegel, Deleuze\index{Deleuze, Gilles (1925--1995)} and Guattari\index{Guattari, F\'elix (1930--1992)} and others makes Plotnitsky's essay an original contribution to the realm of Riemannian philosophy.  Understanding the difference between a philosophical and a mathematical concept is at the center of this essay, like in the previous essay by Jedrezejewski (Chapter 10). Riemann's habilitation lecture,\index{Riemann! habilitation lecture}\index{Habilitation lecture!Riemann} \emph{On the hypotheses that lie at the bases of geometry}, in which mathematics, physics and philosophy are merged, is in the background and provides Plotnitsky with the main material for his argumentation. The question of whether Riemann's notion of space\index{space} belongs to mathematics or to philosophy is central. A notion like the ``plane of immanence" (\emph{plan d'immanence}) as a plane of the movement of thought, in Riemann's approach, is characterized by its multi-component factor, and it is one of the main ways in which Plotnitsky approaches Riemann's work. His essay sheds a new light on Riemann's dissertation and in particular his rethinking of geometry in terms of manifoldness.\index{manifoldness} Connections with works of several philosophers, artists and scientists are highlighted in this chapter. The themes discussed include Leibniz's\index{Leibniz, Gottfried Wilhelm (1646--1716)}  monads,\index{monad} Grothendieck's topoi,\index{topos} and quantum physics.\index{quantum physics} 

\part*{Part III}

Part III of this volume, consisting of Chapters 12 to 16, is mathematical. It covers recent developments in mathematics that are closely related to ideas of Riemann.

Chapter 12, by Feng Luo, is a variation on the Riemann mapping\index{Riemann mapping theorem}\index{theorem!Riemann mapping} theorem and, its generalization, the uniformization theorem. More precisely, it concerns the discrete version of these theorems. 

The interest in a discrete version of the Riemann mapping theorem\index{Discrete Riemann mapping theorem}\index{theorem!discrete Riemann mapping} was given a strong impetus by W. P. Thurston who, in the 1980s, advertised this subject in several lectures and made the relation with circle packings.\index{circle packing}  The idea behind this relation is that a conformal mapping (like the Riemann mapping)\index{Riemann mapping theorem}\index{theorem!Riemann mapping} is characterized by the fact that it sends infinitesimal circles to infinitesimal circles. Circle packings\index{circle packing} involve smaller and smaller circles, therefore they should give information on conformal mappings. An  idea that emerged was that studying circle packings\index{circle packing} might give a new point of view on the Riemann mapping theorem,\index{Riemann mapping theorem}\index{theorem!Riemann mapping} or even a new proof of it. In this setting, a precise question concerning the convergence of circle packings\index{circle packing} to the Riemann mapping theorem was raised and was eventually solved by Rodin and Sullivan.

In Chapter 12, after a presentation of Thurston's ideas on a circle packing version of the Riemann mapping theorem,\index{Riemann mapping theorem}\index{theorem!Riemann mapping} Luo reviews his own recent work on the discrete uniformization theorem\index{uniformization theorem}\index{theorem!uniformization}\index{discrete uniformization theorem}\index{theorem!discrete uniformization} for polyhedral surfaces. The proof is variational. The author highlights relations with approximation theory and with algorithmic and digitalization techniques. 

The material discussed in Chapter 12 may be considered as an illustration of Riemann's ideas on the relation between the discrete\index{discrete} and the continuous,\index{continuous} one of the major themes in his habilitation lecture.\index{Riemann! habilitation lecture}\index{Habilitation lecture!Riemann}
 
The next chapter concerns the Riemann--Roch theorem.
   
The history of the Riemann--Roch Theorem\index{Riemann--Roch theorem}\index{theorem!Riemann--Roch} starts with the so-called Riemann existence theorem, which asserts the existence of meromorphic functions on Riemann surfaces.
The classical Riemann--Roch\index{Riemann--Roch theorem}\index{theorem!Riemann--Roch} theorem gives more precise information. It concerns the dimension of the space\index{space} of meromorphic functions on a compact surface having poles of (at most) a certain order at some prescribed set of points. The theorem is a formula, expressing this dimension in terms of the genus of the surface, thus establishing a fundamental relation between topological and  analytical notions.

 There are several classical proofs of this theorem, some of them topological, others geometric and there are proofs involving abstract algebra, adapted to the case where the ground field (the field of scalars) is more general than that of the complex numbers. The result has many applications, and there are several versions and generalizations of the Riemann--Roch Theorem\index{Riemann--Roch theorem}\index{theorem!Riemann--Roch}. Brill and Noether, back in 1874, already gave an algebro-geometric version of this theorem,\footnote{A. Brill and M. Noether,  \"Uber die algebraischen Functionen und ihre Anwendung in der Geometrie, Math. Ann.  7 (1874,) No. 2, 269-- 316.} a version which is sometimes called the Riemann-Brill-Noether Theorem and which has vast modern developments.
The Riemann-Roch Theorem\index{Riemann--Roch theorem}\index{theorem!Riemann--Roch} was widely generalized by Hirzebruch in 1953, from Riemann surfaces to the setting of projective singular varieties over complex numbers. The modern version of this theorem is expressed in the setting of schemes. Grothendieck\index{Grothendieck, Alexander (1928--2014)} obtained a very general version of the Riemann--Roch\index{Riemann--Roch theorem}\index{theorem!Riemann--Roch} theorem, formulated in the language of categories and functors, which holds for algebraic varieties defined over arbitrary ground fields. This was one of Grothendieck's major discoveries. 
In Chapter 2, \S 2.8 of his \emph{Récoltes et semailles},\footnote{A. Grothendieck, Récoltes et semailles, Réflexions et témoignage sur un passé de mathématicien, unpublished manuscript, 1985--1986, 929 p.} he writes: ``The year 1957 is the one where I was led to extract the theme `Riemann--Roch'\index{Riemann--Roch theorem}\index{theorem!Riemann--Roch}  (Grothendieck's version) which overnight consecrated me `great star'."\footnote{L'année 1957 est celle où je suis amené à dégager le thème ``{Riemann--Roch}" (version Grothendieck)\index{Riemann--Roch theorem}\index{theorem!Riemann--Roch} -- qui, du jour au lendemain, me consacre ``grande vedette".}
 Grothendieck's version of the Riemann--Roch\index{Riemann--Roch theorem}\index{theorem!Riemann--Roch} theorem was the starting point of topological K-theory. 
In the section called \emph{La vision -- ou douze th\`emes pour une harmonie} (The vision -- or twelve themes for a harmony) of \emph{R\'ecoltes et semailles} (Chapter 2, \S\,2.8), Grothendieck considers the theme he calls the \emph{Riemann--Roch-Grothendieck Yoga} as one of the twelve themes which he describes as his ``great ideas" (\emph{grandes id\'ees}).
 There is also a discrete Grothendieck-Riemann--Roch\index{discrete Riemann--Roch theorem}\index{theorem!discrete Riemann--Roch} theorem.\index{Grothendieck-Riemann--Roch theorem}\index{theorem!Grothendieck-Riemann--Roch}
 The famous Atiyah--Singer index theorem,\index{Atiyah--Singer theorem}\index{theorem!Atiyah--Singer}\index{index theorem}\index{theorem!index}  discovered in 1963, can be considered as another generalization of the Riemann--Roch theorem.\index{Riemann--Roch theorem}\index{theorem!Riemann--Roch}

All this justifies the inclusion of a chapter on the Riemann--Roch theorem, whose original idea started with Riemann.

Thus, in Chapter 13, Norbert A'Campo, Vincent Alberge and Elena Frenkel present a modern version of the Riemann--Roch theorem.\index{Riemann--Roch theorem}\index{theorem!Riemann--Roch}  It concerns the space of holomorphic line bundles over a Riemann surfaces. The proof uses Dolbeault cohomology, Serre duality for line bundles, and functional analysis (Fredholm operator theory). The chapter is intended to be a self-contained proof of this cohomological version of Riemann--Roch.\index{Riemann--Roch theorem}\index{theorem!Riemann--Roch} All the required notions (holomorphic line bundle, degree, the Poincaré-Hopf index formula,  the Picard group, sheaves, sheaf cohomology, Chern class, the Cauchy-Riemann operator, Dolbeault cohomology, Serre duality, the index of a Fredholm operator and divisor) are introduced and  clearly explained, in a concise but sufficiently detailed manner so that the reader can understand the theorem and its proof.

The modern version of the Riemann--Roch theorem is an important monument of twentieth-century mathematics.

 Chapter 14, by Victor Pambuccian, Horst Struve and Rolf Struve, concerns the foundations of mathematics. The reader might wonder about the existence of a relation between Riemann and the foundations of mathematics. This relation is hinted on by Riemann in his 1854 habilitation lecture.\index{Riemann! habilitation lecture}\index{Habilitation lecture!Riemann}
 At the beginning of this lecture, Riemann mentions the axiomatic approach as one of the possible approaches to geometry (the other one, to which he will stick soon after, being the metrical approach). He does not further develop this idea, but he raises the issue of the necessity having a \emph{solid foundation} of geometry. Riemann's immediate successors knew about his interest in axiomatics.  W. A. Clifford,\index{Clifford, William Kingdon (1845--1879)} one of the earliest commentators of Riemann's works, writes, in a text titled \emph{The postulates of space},  p. 565\footnote{\label{f:New} Cf. The World of Mathematics,  edited by J. R. Newman, Volume 1, Simon and Schuster, 1956, New York,  552--557.}: ``It was Riemann, however, who first accomplished the task of analysing all the assumptions of geometry, and showing which of them were independent." Helmholtz, in his lecture \emph{On the origin and significance of of geometrical axioms}\footnote{\emph{Ibid.} p. 647--668.} mentions several times Riemann's ideas on the axiomatic foundation of mathematics. 
 It is also useful to recall that Hilbert, in an appendix of his \emph{Foundations of Geometry}, mentions Riemann. He writes: ``The investigations by Riemann and Helmholtz for the foundations of geometry led Lie to take up the problem of the \emph{axiomatic}\footnote{The emphasis is in the original text.} treatment of geometry as introductory to the study of groups." Although Riemann did not develop the axiomatic point of view in any of his own writings, his heirs did, and in particular there were several attempts to axiomatize Riemannian geometry. This is the subject of Chapter 14 of the present volume.

As a matter of fact, the question of the foundations of geometry, like many other foundational questions, can be traced back to Aristotle, developed in his \emph{Posterior analytics} and his other essays. Geometry, at that time, meant mostly Euclidean geometry, although spherical geometry was also known. In any case, the question of the foundations of Riemannian geometry naturally stems from Riemann's work.  With this idea in mind, in Chapter 14, the authors present a set of approaches to the axiomatization of metric spaces, developed by several authors, some of them motivated by Riemann's work. These authors used new notions from various fields that were developed in the few decades that followed Riemann's work: transformation groups, Lie groups, the foundations of arithmetics, mathematical logic, and metric geometry. Although the abstract notion of group is absent from Riemann's writings, the ideas of homogeneity and symmetry are present at several places in his work. The discussion involving group theory that is done in the chapter by Pambuccian, H. Struve and R. Struve is welcome as an important element in the development of Riemann's ideas.

In Chapter 15, Toshikazu Sunada surveys some of the impact of the idea of a Riemann sum --  the basic element of Riemann's integration theory --  in various branches of mathematics. 
He reviews in particular how Riemann sums are used in some counting problems in elementary number theory and in the theory of quasicrystals. The chapter contains illuminating examples, and the author makes interesting connections between works of Riemann Fermat, Dirichlet, Gauss, Siegel, Delone and others.

In Chapter 16, which is the last chapter of Part III, Jacques Franchi gives an exposition of the extension of the theory of Brownian motion\index{Brownian motion}\index{Brownian motion!Riemannnian} to the setting of Riemannian manifolds and of  recent work on relativistic Brownian motion.\index{Brownian motion!relativistic}\index{relativistic Brownian motion}

 We recall that the concept of Brownian motion was introduced initially as a description of the (random) motion of a particle subject to the action of a multitude of other particles in a fluid. Einstein published in 1905 a paper on this subject, in the setting of his kinetic theory of gases. A rigorous mathematical theory of Brownian motion\index{Brownian motion} was developed later, in particular by N. Wiener,\index{Wiener, Norbert (1894--1964)} around the 1920s, on a probabilistic basis and in terms of stochastic processes. 
 We note incidentally that Brownian motion is closely related to the theory of Riemann surfaces. In particular, the Riemann mapping theorem\index{Riemann mapping theorem}\index{theorem!Riemann mapping} can be proved using Brownian motion.\index{Brownian motion} Such an approach was promoted by Sullivan and Thurston. One can also mention, in the same vein, a probabilistic proof of the Riemann mapping theorem by Patodi (1970) and two other proofs by Bismut (1984 and 1985)  of the Atiyah-Singer theorem. These proofs are simpler than the original, using only the Gauss--Bonnet theorem.\index{Gauss--Bonnet theorem}\index{theorem!Gauss-Bonnet} There is also a probabilistic proof of Picard's small theorem by B. Davis (1975).\footnote{I owe the last two examples to Jacques Franchi.} These are only a few of the instances where probability is used to prove results in Riemannian geometry. The theory of Brownian motion\index{Brownian motion} on Riemannian manifolds\index{Riemannnian Brownian motion}\index{Brownian motion!Riemannnian} was developed around 1970 by probabilists. This topic, which is the subject of Chapter 16,  is another occasion for understanding the strong relation between Riemannian geometry and  probability theory.
 
 After his exposition of Brownian motion in Riemannian geometry,  Franchi moves on to the extension of Brownian motion\index{Brownian motion} to a relativistic framework.\index{relativistic Brownian motion}\index{Brownian motion!relativistic} This makes a new relation between Riemannian geometry and relativity theory, and it adds an element  to explain Einstein's strong interest  in this field. As Franchi explains, the relativistic extension of Brownian motion\index{relativistic Brownian motion}\index{Brownian motion!relativistic} is a non-trivial theory, especially because of the relativistic constraint that the particle's velocity cannot exceed that of light. 

 Chapter 16 is also the occasion of following the history of the interesting theory of diffusion, where the first (negative) results were obtained by Dudley in 1965, who proved that a Lorentz-covariant Markov diffusion process cannot exist in the framework of special relativity, in particular because of the same problem of large velocities. At the same time, Dudley proposed a construction of a relativistic diffusion at the level of the tangent bundle of Minkowski space. He specified the asymptotic behavior of that diffusion and he showed that it is canonical, given the constraints of being covariant under the action of the Lorentz group.
A similar approach on the unit tangent bundle of a generic Lorentzian manifold, that is, in the setting of general relativity,\index{general relativity}  was made by Franchi and Le Jan in 2007. In this setting, relativistic diffusion\index{relativistic diffusion} becomes a random perturbation of the geodesic flow over a Lorentzian manifold.\index{Lorentzian manifold} Some basic examples are then analyzed to some extent.

 The exposition in Chapter 16 follows the gradual move from the Euclidean to the Riemannian and then to the relativistic worlds. This theory is another instance of the intricate interaction between geometry, analysis, probability and physics relying heavily on Riemann's ideas. Thus, this chapter makes a natural transition between Part III and  the next part of the book.

\part*{Part IV}

 Part IV of this volume concerns physics. It contains three chapters on the extension of Riemann's ideas to modern physics, mainly, to relativity theory. Riemann, in his habilitation lecture\index{Riemann! habilitation lecture}\index{Habilitation lecture!Riemann}  \emph{\"Uber die Hypothesen, welche der Geometrie zu Grunde liegen}, expressed the fact that physical space might not satisfy the axioms of Euclidean geometry. This, is, from the philosophical point of view, the starting point of Riemann's position as a predecessor of modern physics. At a more practical level, the mathematical development of Einstein's theory of  general relativity is, at a fundamental level, in the tradition of Riemann's differential geometry; this is one of the themes of Part IV of the present volume.
 
Minkowski geometry is a semi-Riemannian geometry where the metric tensor is not positive-definite -- a mathematical consequence of the physical fact that particles cannot move at a speed larger than that of light. Although the geometric setting of special relativity is Minkowski geometry,\index{Minkowski geometry} which is not Riemannian, the basic mathematical ideas that are used in the development of this geometry are similar to those introduced by Riemann.
In other words, the fact that Minkowski geometry  differs from Riemannian geometry does not affect the fact that it is in the lineage of Riemann's ideas on geometry and space.\index{space}   Riemann's discussion of the invariance properties of a metric, which he carries in his habilitation lecture, have their analogues as invariance properties of the Lorentz transformations of special relativity. In fact, many important features of a four-manifold equipped with a Riemannian metric together with its Riemann's curvature tensor have their analogue in Minkowski spacetime.\index{Minkowski spacetime} In general relativity, the metric tensor that describes the local geometry of space is the mathematical representation of the gravitational potential. This is again in the tradition of Riemann, who conceived his infinitesimal metric tensor in close relation with physics, the curvature of the space being seen as a consequence of the physical infinitesimal forces.

Andreas Hermann and Emmanuel Humbert, in Chapter 17 of this volume, study a variant of the so-called Positive Mass Conjecture\index{conjecture!positive mass}\index{positive mass conjecture} for closed Riemannian manifolds.
The conjecture is a statement in general relativity which gives conditions under which  the mass of an asymptotically flat spacetime is non-negative. The relevance of this theory to the subject of this book is that it is an important instance of a purely physical problem that can be formulated in terms of Riemannian geometry.  The subject discussed has a long history. Using minimal surfaces and variational methods, Schoen and Yau proved in 1979 the positive mass conjecture for 3-dimensional Riemannian manifolds, in the setting of the Hamiltonian formulation of general relativity. Witten gave a later proof (1981) which holds in any dimension. Roughly speaking, the positive mass theorem says that if the scalar curvature of a spacetime is everywhere positive, then its mass is positive. An inequality attributed to Penrose says that the mass of a
spacetime can be estimated by the total area of the black holes contained in it, and that equality is attained only for a simple model of a black hole, the so-called Schwarzschild model.\index{Schwarzschild model}

 Chapter 18, by Marc Mars,  focuses on an important aspect of  the rich interaction between mathematics and physics based on the interplay between differential geometry, in the tradition of Riemann, and gravity, in the setting of the theory of relativity.  This aspect is the
local characterization of pseudo-Riemannian manifolds, which is central in general relativity in order to identify spacetime geometries independently of the
specific set of coordinates used to describe them. One of the many groundbreaking contributions of Riemann to geometry was the introduction of his tensor (the so-called Riemann tensor) which vanishes if and only if the metric is locally flat. It turns out that this fundamental local characterization result holds
independently of the signature of the metric, and is
the motivation of many other similar characterization theorems.

 After reviewing the classical results on the subject, including 
\'Elie Cartan's\index{Cartan, \'Elie Joseph (1869--1951)} characterization of Riemannian locally symmetric spaces in terms of the 
derivative of  the Riemann tensor
 and  Weyl's\index{Weyl, Hermann (1885--1955)} characterization of locally conformally flat spaces in terms of the vanishing of the
conformal curvature (i.e.,
Weyl's) tensor\index{Weyl tensor}, the chapter discusses a selection of various characterization results of physically relevant spacetimes.
The emphasis is primarily on spacetimes
describing stationary black holes, both in the static and in the rotating case. Thus, several characterization results are presented for the
Schwarzschild\index{Schwarzschild spacetime} (and Kruskal) 
spacetimes, as well as for the Kerr metric, and its charged and cosmological constant generalizations.
Local characterization of other spacetimes, such as for instance $pp$-waves and related spacetimes, are also described.

 The last chapter of Part IV, Chapter 19, by Jean-Philippe Nicolas, contains an exposition of Penrose's conformal technique and its application to asymptotic analysis\index{asymptotic analysis} in general relativity. The setting is again that of Lorentzian geometry\index{Lorentzian geometry}\index{geometry!Lorentzian} considered as an extension of Riemannian geometry in which space and time are united by an indefinite metric of signature (1, 3). The author presents Penrose's approach to general relativity with the central role played by the light cone structure and he explains its relation with Riemannian geometry and with Einstein's theory. The focus is on Penrose's use of conformal compactifications\index{conformal compactification} in the study he made of the asymptotic properties of spacetimes and fields. Indeed, Penrose introduced in 1963 a basic geometrical construction which is termed in Chapter 19 a ``compactified unphysical" spacetime. This is a manifold with boundary to which the conformal metric extends smoothly, i.e., there is a metric in the conformal class that extends as a smooth non-degenerate Lorentzian metric. Spacetimes that admit smooth conformal compactifications are characterized by a decay property of their Weyl curvature\index{Weyl curvature}\index{curvature!Weyl} at infinity. When such a compactification exists, the boundary of the manifold is equipped with a nice geometric stratified structure.
 
After surveying Penrose's theory, Nicolas reviews some of its applications to questions of scattering\index{scattering} and peeling.\index{peeling} Scattering theory is a way of studying the evolution of solutions of a certain equation by a so-called scattering operator, an operator which associates to the asymptotic behavior of the solutions in the distant past their asymptotic behavior in the distant future. Peeling is a generic asymptotic behavior discovered by R. Sachs in the beginning of the 1960s. In the mid 1960s, Penrose proved that this behavior is equivalent to the boundedness of the rescaled field at infinity, using the conformal method and the $2$-spinor formalism. The question of the genericity of the peeling behavior is discussed. A new approach to these questions together with results by L. J. Mason and Nicolas are presented. The two approaches to asymptotic analysis described in Chapter 19 make a fundamental use of the notion of conformal compactification. 

Beyond the results presented, the true focus of the essay is on the nature of spacetime: whereas many modern approaches to general relativity break the symmetry between time and space by performing a $3+1$ splitting of the geometry, Penrose's approach truly deals with the $4$-dimensional manifold and relies on causal objects like lightcones instead of Cauchy hypersurfaces.

  \medskip

  Reading the texts of the ancient mathematicians always sheds a new light on the problems that nurture us every day. Regarding Riemann, Weil writes the following in his \emph{Apprenticeship of a Mathematician}:\footnote{A. Weil, The Apprenticeship of a Mathematician, Springer, Basel, 1991. Translated from the French: Souvenirs d'apprentissage, Basel, Birkh\"auser, 1991, p. 40.}
\begin{quote}\small
[...] In the same year, I began to read Riemann. Some time earlier, and first of all in reading Greek poets, I had become convinced that what really counts in the history of humanity are the truly great minds, and that the only way to get to know these minds was through direct contact with their works. I have since learned to modify this judgement quite a bit, though I have never really let go of it completely. My sister, however, who had come to a similar viewpoint -- either on her own or perhaps partly under my influence -- held on to it until the very end of her too short life. During my year of instruction in philosophy, I had also been struck by a phrase of Poincar\'e's\index{Poincaré, Henri (1854--1912)} which expresses no less an extreme position: ``The value of civilizations lies only in their sciences and arts."  With such ideas in my mind I had no choice but to dive headlong into the works of the great mathematicians of the past, as soon as they were materially and intellectually within my grasp. Riemann was the first; I read his inaugural dissertation and his major work on Abelian functions. Starting out thus was a stroke of luck of which I have always been grateful. These are not hard to read as long as one realizes that every word is loaded with meaning: there is perhaps no other mathematician whose writing matches Riemann's for density. Jordan's second volume was good preparation for studying Riemann. Moreover, the library\footnote{The library is meant to be that of the \'Ecole Normale Supérieure.} had a good collection of Felix Klein's mimeographed lecture notes, a large part of which is simply a rather discursive, but intelligent, commentary fleshing out of the extreme concision of Riemann's work.
 \end{quote}

\medskip
  The theme of the present volume, beyond the reference to Riemann's work, belongs to the more general profound interrelation between mathematics, physics and philosophy. The relation is multiple. Physics may exert an influence on mathematics and vice versa. Physics has also an impact on philosophy, and philosophy on mathematics. The ancient Greeks, the founders of mathematics as a deductive science in the way we intend it today, were completely aware of these interrelations. One may mention here Archimedes, Ptolemy and many others great figures.
Euler, Poincaré\index{Poincaré, Henri (1854--1912)} and Cauchy were also physicists and philosophers, and they also wrote on the interrelations and the impact of these fields, each of them with his own style and according to his own interests. The subject of Euler's first public lecture, delivered in Basel in 1724 (the year Immanuel Kant was born), at the occasion of his obtention of his philosophy diploma, was the comparison between the philosophical systems of Newton and Descartes. Euler's philosophy is at some places religiously oriented and some of his philosophical writings are permeated with theological considerations. They were influential to his approach to physics. We allude to this and we give some examples in Chapter 1 of the present volume.  The philosophical writings of Cauchy, who like Euler, was a devote Christian, and who was furthermore involved in several charities, are also infused with religion. In several passages, he mentions the limitation of mathematics. For instance, in the introduction to his \emph{Cours d'analyse},\footnote{A.-L. Cauchy, Cours d'analyse de l'\'Ecole Royale Polytechnique, 1${}^{\mathrm{re}}$ partie. Analyse alg\'ebrique. Imprimerie royale,  Paris, 1821.
\OE uvres compl\`etes, s\'erie 2, tome III.} he writes the following:
\begin{quote}\small
[...] Thus, let us be persuaded that there are other truths than those of algebra, realities other than sensible objects. Let us cultivate ardently the mathematical sciences, without trying to extend them beyond their domain; and let us not imagine that one can tackle history with formulae, or use theorems of algebra or of differential calculus as an assent to morals.\footnote{[...] Soyons donc persuadés qu'il existe des vérités autres que les vérités de l'algèbre, des réalités autres que les objets sensibles. Cultivons avec ardeur les sciences mathématiques, sans vouloir les étendre au-delà de leur domaine ; et n'allons pas nous imaginer qu'on puisse attaquer l'histoire avec des formules, ni donner pour sanction à la morale des théorèmes d'algèbre ou de calcul intégral.}
\end{quote}

Regarding physics and its development, we quote the following passage, from lectures that Cauchy gave in Turin in 1833,  \emph{Sept le\c cons de physique g\'en\'erale}, (Seven lessons on general physics), p. 5:\footnote{A.-L. Cauchy, Sept le\c cons de physique g\'en\'erale, Paris, bureau du journal Les Mondes  et Gauthier-Villars, 1868.}
\begin{quote}\small
[...] Among these sciences, there is one in which all the power of analysis is manifested, and in which calculus, created by man, takes care 
of teaching him, through a mysterious language, the links that exist between phenomena which apparently are very different, and between the particular and the general laws of creation. This science, which we can trace back to the discovery of the principle of universal gravitation, was successively enriched by the immortal works of people like Descartes,\index{Descartes, Ren\'e (1596--1650)}  Huygens,\index{Huygens, Christiann (1629--1695)}  Newton, \index{Newton, Isaac (1643--1727)} and Euler.\index{Euler, Leonhard (1707--1783)}  But it is particularly since twenty years that the rapid improvement of mathematical analysis allowed him to make huge progress. It is since that epoch that we were able to apply calculus to the theory of elasticity, to that of heat propagation in solids or in space, of the propagation of waves on the surface of a heavy fluid, of the transmission of sound through solid bodies;  to the theory of dynamical elasticity, to that of vibration of plates or elastic lamina; and finally to the theory of light including the various reflection phenomena, simple refraction, double refraction, polarization, colors, etc. Finally, it is since that epoch that important works of people like Amp\`ere,\index{Amp\`ere, Andr\'e-Marie (1775--1836)} Fourier,\index{Fourier, Joseph (1768--1830)}  Poisson\index{Poisson, Sim\'eon Denis (1781--1840)}  and of some others of which I do not need  to remind  you the names, were published.\footnote{ [...] Parmi ces sciences, il est une où se manifeste toute la puissance de l'analyse, et dans laquelle le calcul cr\'e\'e par l'homme se charge de lui apprendre, par un myst\'erieux langage, les liaisons qui existent entre des ph\'enom\`enes en apparence tr\`es divers, entre les lois particuli\`eres et les lois g\'en\'erales de la cr\'eation. Cette science, qu'on peut faire monter à la d\'ecouverte du principe de la gravitation universelle, a \'et\'e successivement enrichie des immortels travaux des Descartes,\index{Descartes, Ren\'e (1596--1650)} des Huyghens, des Newton,\index{Newton, Isaac (1643--1727)} des Euler. Mais c'est particuli\`erement depuis vingt ans que le perfectionnement rapide de l'analyse math\'ematique lui a permis de faire d'immenses progr\`es. C'est depuis cette \'epoque qu'on a pu appliquer le calcul à la th\'eorie de l'\'elasticit\'e, de la propagation de la chaleur dans des corps ou dans l'espace, de la propagation des ondes à la surface d'un fluide pesant, de la transmission du son à travers les corps solides ;
  à la th\'eorie de l'\'elasticit\'e dynamique, à celle des vibrations des plaques ou des lames \'elastiques ; enfin à la th\'eorie de la lumi\`ere comprenant les ph\'enom\`enes divers de la r\'eflexion, de la r\'efraction simple, de la double r\'efraction, de la polarisation, de la coloration, etc... C'est enfin depuis cette \'epoque qu'ont \'et\'e publi\'es les importants travaux des Amp\`ere, des Fourier, des Poisson et de quelques autres dont il est inutile de vous rappeler les noms.} 
\end{quote}

To stay close to the epoch of Riemann, we quote another one of his close predecessors, Joseph Fourier,\index{Fourier, Joseph (1768--1830)} from his \emph{Théorie analytique de la chaleur} (Analytic theory of heat), published in 1822, a text which was very important for Riemann who refers to it in his habilitation memoir on trigonometric functions. On the relation between mathematics and the study of nature, Fourier writes in  the Introduction to his work:\footnote{J. Fourier, Théorie analytique de la chaleur, Paris, Firmin Didot, 1822, Discours préliminaire, p. xiii.}
  \begin{quote}\small
  The thorough study of nature is the most profound productive source of mathematical discoveries. Not only this study, offering to the researches a specific purpose, has the advantage of excluding fuzzy questions and dead-end calculations; it is also a secure way of forming the heart of analysis, and of discovering there the elements whose knowledge is the most important to us, and which this science must always preserve: these fundamental elements are those which reproduce themselves in every natural effect.
  One can see, for instance that the same expression, which geometers had considered as an abstract property, and which from this respect belong to general analysis, also represents the motion of light in the atmosphere, that it determines the laws of diffusion of heat in solid matter, and that it enters in the main questions of the theory of probability.\footnote{L'étude approfondie de la nature est la source la plus féconde des découvertes mathématiques. Non seulement cette étude, offrant aux recherches un but déterminé, a l'avantage d'exclure les questions vagues et les calculs sans issue ; elle est encore un moyen assuré de former l'analyse elle-m\^eme, et d'en découvrir les éléments qu'il nous importe le plus de conna\^\i tre, et que cette science doit toujours conserver : ces éléments fondamentaux sont ceux qui se reproduisent dans tous les effets naturels.
  On voit, par exemple, qu'une m\^eme expression, dont les géomètres avaient considéré les propriétés abstraites, et qui sous ce rapport appartient à l'analyse générale, représente aussi le mouvement de la lumière dans l'atmosphère, qu'elle détermine les lois de la diffusion de la chaleur dans la matière solide, et qu'elle entre dans les questions principales de la théorie des probabilités.}
  \end{quote}

Poincaré\index{Poincaré, Henri (1854--1912)} was a prototype of the scientist-philosopher, and it was probably under his influence that most of the pre-eminent French mathematicians of his epoch became deeply interested in physics and philosophy.  We quote him from his 1908 ICM talk (Rome):\footnote{H. Poincaré, l'Avenir des mathématiques,  Atti del IV congresso internazionale dei matematici, Volume 1, Accademia dei Lincei, Rome, 1909,  p. 167--182.}
\begin{quote}\small
We cannot forget what our goal should be. As I see it, it is twofold. Our science borders at the same time on philosophy and on physics, and it is for our two neighbors that we are working. On the other hand, we have always seen, and we shall also see the mathematicians walking in two opposite directions. On the one hand, mathematical science must reflect on itself, and this is useful, because reflecting on itself means reflecting on the human mind that created it, all the more since this is his creation for which he borrowed the less from outside. 
This is why certain mathematical speculations are useful, like the ones which aim at the study of postulates, of unusual geometries, of functions with strange behavior. The more these speculations deviate from the most common conceptions, and consequently of the nature of their applications, the better they will show what human mind is able to do, when it avoids more and more  the tyranny of the external world, and consequently, the more they will let us know it itself. But on the other hand, it is on the side of nature that we must direct the greater part of our army.\footnote{Nous  ne  pouvons oublier 
quel  doit  \^etre  notre  but ;  selon  moi  ce  but  est  double ; notre  science  confine  à la  fois  
à  la  philosophie  et  à  la  physique,  et  c'est  pour  nos  deux voisines  que  nous travaillons ; 
aussi  nous  avons  toujours  vu  et  nous verrons  encore les  mathématiciens  marcher  dans  
deux  directions  opposées.  
D'une  part,  la  science   mathématique  doit   réfléchir  sur  elle-m\^eme   et   cela   est   
utile,
  parce  que  réfléchir  sur  elle-m\^eme,  c'est  réfléchir  sur  l'esprit   humain   qui   l'a   
créée,  d'autant  plus  que  c'est  celle  de  ses  créations  pour laquelle  il  a  fait  le  moins  
d'emprunts  au  dehors.  
C'est  pourquoi  certaines  spéculations  mathématiques  sont utiles, 
comme  celles   qui   visent   l'étude   des   postulats,  des
  géométries
  inaccoutumées,  des  
fonctions  à  allures  étranges. 
 Plus  ces spéculations s'écarteront  des  conceptions  les plus 
communes,  et  par  conséquent   de   la   nature   et   des   applications,  mieux  elles  nous  
montreront  ce  que  l'esprit  humain  peut  faire,  quand  il  se  soustrait  de  plus  en  plus  
à  la  tyrannie  du  monde  extérieur,  mieux  par  conséquent  elles  nous le  feront  connaître  
en lui-m\^eme.  
Mais  c'est  du  côté  opposé,  du  côté  de  la  nature,  qu'il   faut   diriger  le  gros  de  
notre  armée.} 
\end{quote}
 
Closer to us, Grothendieck,\footnote{The fact that Grothendieck\index{Grothendieck, Alexander (1928--2014)} was not interested in physics is a myth. It suffices to read his non-mathematical writings to be convinced of the contrary.} who is quoted several times in the present volume and who at several places declared that he was a heir of Riemann, has a huge amount of still unpublished philosophical writings. In his \emph{Récoltes et semailles}, which we already mentioned in this introduction, expressing his ideas about a ``unitary theory" in physics, and after a long digression involving Euclid, Newton, Riemann and Einstein, Grothendieck writes (\S 2.20, Note 71):  
\begin{quote}\small

To summarize, I foresee that the long-awaited  renewal (if ever it comes...) will rather come from someone who has the soul of a mathematician, who is well informed about the great problems of physics, rather than from a physicist. But above all, we need a man having the ``philosophical openness" that is required to grasp the crux of the problem. The latter is not at all of a technical nature, but it is a fundamental problem of ``natural philosophy."\footnote{Pour r\'esumer, je pr\'evois que le renouvellement attendu (s'il doit encore venir \ldots) viendra plutôt d'un math\'ematicien dans l'\^ame, bien inform\'e des grands probl\`emes de la physique, que d'un physicien. Mais surtout, il y faudra un homme ayant ``l'ouverture philosophique" pour saisir le n\oe ud du probl\`eme. Celui-ci n'est nullement de nature technique, mais bien un probl\`eme fondamental de ``philosophie de la nature."}
\end{quote}

\medskip

  The present volume is a modest tribute to all those who taught us creative science.
  
   \bigskip

\end{document}